\documentclass[11pt]{article}
\begin{document}
\vspace*{.5cm}
\begin{center}
{\Large{\bf   Slant Riemannian maps to K\"{a}hler manifolds}}\\
\vspace{.5cm}
 { Bayram \d{S}ahin} \\
\end{center}

\vspace{.5cm}
\begin{center}
{\it Inonu University, Department of Mathematics, 44280,
Malatya-Turkey. E-mail:bayram.sahin@inonu.edu.tr}
\end{center}
\vspace{.5cm}

\noindent {\bf Abstract.} {\small We introduce slant Riemannian maps
from Riemannian manifolds to almost Hermitian manifolds as a
generalization of slant immersions, invariant Riemannian maps and
anti-invariant Riemannian maps. We give examples, obtain
characterizations and investigate the harmonicity of such maps. We
also obtain necessary and sufficient conditions for slant Riemannian
maps to be totally geodesic. Moreover we relate the notion of slant
Riemannian maps to the notion of pseudo horizontally
weakly conformal maps which are useful for proving various complex-analytic properties of stable harmonic maps from complex projective space. }\\



\pagestyle{myheadings}

\section*{1.~Introduction}

  \setcounter{equation}{0}
\renewcommand{\theequation}{1.\arabic{equation}}
\markboth{Slant Riemannian maps}{Slant Riemannian
maps}{\thispagestyle{plain}}

Smooth maps between Riemannian manifolds are useful for comparing
geometric structures between two manifolds. Isometric immersions
(Riemannian submanifolds) are basic such maps between Riemannian
manifolds and they are characterized by their Riemannian metrics and
Jacobian matrices. More precisely, a smooth map $F:(M_1,
g_1)\longrightarrow (M_2, g_2)$ between Riemannian manifolds $(M_1,
g_1)$ and $(M_2, g_2)$ is called an isometric immersion if $F_*$ is
injective and
\begin{equation}
g_2(F_*X, F_*Y)=g_1(X, Y)\label{eq:1.1}
\end{equation} for $X, Y $
vector fields tangent to $M_1$, here $F_*$ denotes the derivative
map.\\

Let $\bar{M}$ be a K\"{a}hler manifold with complex structure $J$
and $M$ a Riemannian manifold isometrically immersed in $\bar{M}$.
We note that many types of submanifolds can be defined depending on the
behaviour of the tangent bundle of the submanifold under the
action of the complex structure of the ambient manifold. A
submanifold $M$ is called holomorphic (complex) if $J (T_p M)\subset
T_p M$, for every $p \in M$, where $T_p M$  denotes the tangent
space to $M$ at the point $p $. $M$ is called totally real if $J(T_p
M) \subset T_p M^{\perp}$ for every $p \in M,$ where $T_p M^{\perp}$
denotes the normal space to $M$ at the point $p$.   As a
generalization of holomorphic and totally real submanifolds, slant
submanifolds were introduced by Chen in \cite{Chen}. We recall that
the submanifold $M$ is called slant if for all
    non-zero vector $X$ tangent to $M$ the angle $\theta (X)$
    between $J X$ and $T_p M$ is a constant, i.e, it does not depend
    on the choice of $p \in M$ and $X \in T_p M$.\\

On the other hand, Riemannian submersions between Riemannian
manifolds were studied by O'Neill \cite{O'Neill} and Gray
\cite{Gray}. Later such submersions were considered between
manifolds with differentiable structures. As an  analogue of
holomorphic submanifolds, Watson defined almost Hermitian
submersions between almost Hermitian manifolds and he showed that
the base manifold and each fiber have the same kind of structure as
the total space, in most cases \cite{Watson}. We note that almost
Hermitian submersions have been extended to the almost contact
manifolds \cite{Domingo}, locally conformal K\"{a}hler
manifolds\cite{lck}
and quaternion K\"{a}hler manifolds \cite{Ianus} (see:\cite{Falcitelli} for details cocerning Riemannian submersions between Riemannian manifolds equipped with additional structures of complex, contact, locally conformal or quaternion K\"{a}hler). \\

In 1992, Fischer  introduced Riemannian maps between Riemannian
manifolds in \cite{Fischer} as a generalization of the notions of
isometric immersions and Riemannian submersions. Let $F:(M_1,
g_1)\longrightarrow (M_2, g_2)$  be a smooth map between Riemannian
manifolds such that $0<rank F<min\{ m, n\}$, where $dimM_1=m$ and
$dimM_2=n$. Then we denote the kernel space of $F_*$ by $kerF_*$ and
consider the orthogonal complementary space
$\mathcal{H}=(kerF_*)^\perp$ to $kerF_*$. Then the tangent bundle of
$M_1$ has the following decomposition

$$TM_1=kerF_* \oplus\mathcal{H}.$$

We denote the range of $F_*$ by $rangeF_*$ and consider the
orthogonal complementary space $(rangeF_*)^\perp$ to $rangeF_*$ in
the tangent bundle $TM_2$ of $M_2$. Since $rankF<min\{ m, n\}$, we
always have $(rangeF_*)^\perp$. Thus the tangent bundle $TM_2$ of
$M_2$ has the following decomposition
$$TM_2=(rangeF_*)\oplus (rangeF_*)^\perp.$$
 Now, a smooth
map $F:(M^{^m}_1,g_1)\longrightarrow (M^{^n}_2, g_2)$ is called
a Riemannian map at $p_1 \in M$ if the horizontal restriction
$F^{^h}_{*p_1}: (ker F_{*p_1})^\perp \longrightarrow (range
F_{*p_1})$  is a linear isometry between the inner product spaces
$((ker F_{*p_1})^\perp, g_1(p_1)\mid_{(ker F_{*p_1})^\perp})$ and
$(range F_{*p_1}, g_2(p_2)\mid_{(range F_{*p_1})})$, $p_2=F(p_1)$.
Therefore Fischer stated in \cite{Fischer} that a Riemannian map is
a map which is as isometric as it can be. In another words, $F_*$
satisfies the equation (\ref{eq:1.1}) for $X, Y$ vector fields
tangent to $\mathcal{H}$. It follows that isometric immersions and
Riemannian submersions are particular Riemannian maps with
$kerF_*=\{ 0 \}$ and $(rangeF_*)^\perp=\{ 0 \}$. It is known that a
Riemannian map is a subimmersion \cite{Fischer}.  It is also important to note that Riemannian maps satisfy the eikonal equation which is a bridge between geometric optics and physical optics. For Riemannian maps and their applications in spacetime geometry see: \cite{Garcia-Kupeli}.\\

In \cite{Sahin1}, we introduced invariant and anti-invariant
Riemannian maps as a generalization of holomorphic submanifolds and
totally real submanifolds of almost Hermitian manifolds. We gave
examples of such maps and obtained fundamental properties of
invariant and anti-invariant Riemannian maps. As a generalization of invariant and anti-invariant Riemannian maps, semi-invariant Riemannian maps from Riemannian manifolds to almost Hermitian manifolds were defined and the geometry of such maps was studied in \cite{Sahin2}.\\

In this paper, we introduce slant Riemannian maps from Riemannian
manifolds  to almost Hermitian manifolds. We show that slant
Riemannian maps include slant immersions (therefore holomorphic
immersions and totally real immersions), invariant Riemannian maps
and anti-invariant Riemannian maps. We also obtain an example which
is not included in immersions, or invariant Riemannian maps or
anti-invariant Riemannian maps. We investigate the harmonicity of
slant Riemannian maps and obtain necessary and sufficient conditions
for such maps to be totally geodesic. We also show that every slant
Riemannian map is a pseudo horizontally weakly conformal map, then
we obtain necessary and sufficient conditions for slant Riemannian
maps to be a pseudo homothetic map. We note that the notion of
pseudo horizontally weakly conformal maps were introduced in
\cite{BBBR} to study the stability of harmonic maps into
irreducible Hermitian symmetric spaces of compact type, later such
maps have been studied in \cite{AA}, \cite{AAB} and \cite{LS}.

\section*{\bf 2.~Preliminaries}
\setcounter{equation}{0}
\renewcommand{\theequation}{2.\arabic{equation}}

In this section, we recall basic facts about Riemannian maps and
almost Hermitian manifolds. Let $(M, g_{_M})$ and $(N, g_{_N})$ be
Riemannian manifolds and suppose that $F: M\longrightarrow N$ is a
smooth map between them. Then the differential $F_*$ of $F$  can be viewed as a section of the bundle $Hom(TM, F^{-1}TN)\longrightarrow M,$
where $F^{-1}TN$ is the pullback bundle which has fibres
$(F^{-1}TN)_p=T_{F(p)} N, p \in M.$ $Hom(TM, F^{-1}TN)$ has a
connection $\nabla$ induced from the Levi-Civita connection
$\nabla^M$ and the pullback connection. Then the second fundamental
form of $F$ is given by
\begin{equation}
(\nabla {F}_*)(X, Y)=\nabla^{F}_X {F}_*(Y)-{F}_*(\nabla^M_X Y)
\label{eq:2.1}
\end{equation}
for $X, Y \in \Gamma(TM)$, where $\nabla^F$ is the pullback
connection. It is known that the second fundamental form is
symmetric. Let $F: (M, g_{_M}) \longrightarrow (N, g_{_N})$ be a
smooth map between Riemannian manifolds and assume that $M$ is compact,
then its energy is
\begin{eqnarray}
E(F)=\int_{M}e(F)\,v_{g}=\frac{1}{2}\int_{M}|dF|^{2}\,v_{g}.\nonumber
\end{eqnarray}
The critical points of $E$  are called harmonic maps. Standard
arguments yield the associated Euler-Lagrange equation, the
vanishing of the tension field: $\tau(F)=trace (\nabla
{F}_*)$ \cite{Baird-Wood}. \\

From now on, for simplicity, we denote by $\nabla^2$ both the
Levi-Civita connection of $(M_2, g_2)$ and its pullback along $F$.
 Then according to \cite{Nore}, for any vector field $X$ on $M_1$ and any section $V$ of $(range F_*)^\perp$, where $(range F_*)^\perp$ is the
 subbundle of $F^{-1}(TM_2)$ with fiber $(F_*(T_pM))^\perp$-orthogonal complement of $F_*(T_pM)$ for $g_2$ over $p$, we have
 $\nabla^{^{F \perp}}_XV$ which is the orthogonal projection of $\nabla^2_XV$ on $(F_*(TM))^\perp$. In \cite{Nore}, the author also showed that $\nabla^{^{F \perp}}$ is a linear connection on $(F_*(TM))^\perp$ such that $\nabla^{^{F \perp}}g_2=0$. We now  define $\mathcal{S}_{V}$ as
\begin{equation}
\nabla^2_{_{F_*X}}V=-\mathcal{S}_{_V}F_*X+\nabla^{^{F
\perp}}_{_{X}}V, \label{eq:2.2}
\end{equation}
where $\mathcal{S}_{_V}F_*X$ is the tangential component (a vector
field along $F$) of $\nabla^2_{_{F_*X}}V$. It is easy to see that
$\mathcal{S}_V F_*X$ is bilinear in $V$ and $F_*X$ and
$\mathcal{S}_V F_*X$ at $p$ depends only on $V_p$ and $F_{*p}X_p$.
By direct computations, we obtain
\begin{equation}
g_2(\mathcal{S}_{_V} F_*X,F_*Y)=g_2(V, (\nabla F_*)(X,Y)),
\label{eq:2.3}
\end{equation}
for $X, Y \in \Gamma((ker F_*)^\perp)$ and $V \in \Gamma((range
F_*)^\perp)$. Since $(\nabla F_*)$ is symmetric, it follows that
$\mathcal{S}_{_V}$ is a symmetric linear transformation of $range F_*$.\\

A $2k$-dimensional Riemannian manifold $(\bar{M}, \bar{g}, \bar{J})$
is called an almost Hermitian manifold if there exists a tensor
filed  $\bar{J}$ of type (1,1) on $\bar{M}$ such that $\bar{J}^2 = -
I$ and
\begin{equation}
\bar{g}(X,Y)=\bar{g}(\bar{J}X,\bar{J}Y) \label{eq:2.4}, \forall X,Y
\in \Gamma(T\bar{M}),
\end{equation}
 where $I$ denotes the identity transformation of
$T_{p} \bar{M}$. Consider an almost Hermitian manifold $(\bar{M},
\bar{J}, \bar{g})$  and denote by $\bar{\nabla}$ the Levi-Civita
connection on $\bar{M}$ with respect to $\bar{g}$, then $\bar{M}$ is
called a K\"{a}hler manifold \cite{Yano-Kon} if $\bar{J}$  is
parallel with respect to $\bar{\nabla}$, i.e,
\begin{equation}
(\bar{\nabla}_X \bar{J})Y=0 \label{eq:2.5}
\end{equation}
for $X, Y \in \Gamma(T\bar{M}).$

\section*{3. Slant Riemannian maps}
  \setcounter{equation}{0}
\renewcommand{\theequation}{3.\arabic{equation}}

In this section, we define slant Riemannian maps from a Riemannian
manifold to an almost Hermitian manifold. We give examples, obtain
characterizations and  investigate the harmonicity of a slant
Riemannian map. Then we give necessary and sufficient conditions
for a slant Riemannian map to be totally geodesic. \\

\noindent{\bf Definition~3.1.~}{\it Let $F$ be a Riemannian map from
a Riemannian manifold $(M_1,g_1)$ to an almost Hermitian manifold
$(M_2,g_2,J)$. If for any non-zero vector $X \in \Gamma(ker
F_*)^\perp$, the angle $\theta(X)$ between $JF_*(X)$ and the space
$range F_*$ is a constant, i.e. it is independent of the choice of
the point $p \in M_1$ and choice of the tangent vector $F_*(X)$ in
$range F_*$,
then we say that $F$ is a slant Riemannian map. In this case, the angle $\theta$ is called the slant angle of the slant Riemannian map.}\\

Since $F$ is a subimmersion, it follows that the rank of $F$ is
constant on $M_1$, then the rank theorem for functions implies that
$ker F_*$ is an integrable subbundle of $TM_1$, (\cite{Marsden},
page:205). \\

We first give some examples of slant Riemannian maps.\\

\noindent{\bf Example~1.~}Every $\theta-$ slant immersion from a
Riemannian manifold to an almost Hermitian manifold is a slant
Riemannian map with  $ker F_*=\{0\}$.\\

\noindent{\bf Example~2.~}Every invariant Riemannian map from a
Riemannian manifold to  an almost Hermitian manifold is a slant Riemannian map with $\theta=0$.\\

\noindent{\bf Example~3.~}Every anti-invariant Riemannian map from a
Riemannian manifold to  an almost Hermitian manifold is a slant Riemannian map with $\theta=\frac{\pi}{2}$.\\

A slant Riemannian map is said to be proper if it is not an
immersion and $\theta\neq 0, \frac{\pi}{2}$. Here is an example of proper slant Riemannian maps.\\

 \noindent{\bf Example~4.~}Consider the following Riemannian
map given by
$$
\begin{array}{cccc}
  F: & R^4             & \longrightarrow & R^4\\
     & (x_1,x_2,x_3,x_4) &             & (x_1,\frac{x_2+x_3}{\sqrt{3}},\frac{x_2+x_3}{\sqrt{6}},0).
\end{array}
$$
Then $rank F=2$ and for any $0< \alpha <\frac{\pi}{2}$, $F$ is a
slant Riemannian map
with slant angle $\cos^{-1}(\sqrt{\frac{2}{3}})$.\\

We also have the following result which is based on the fact that
the composition of a Riemannian submersion $F_1$ from a Riemannian
manifold $(M_1,g_1)$ onto a Riemannian manifold $(M_2,g_2)$ and an
isometric immersion $F_2$ from $(M_2,g_2)$ to a Riemannian manifold
$(M_3,g_3)$ is a Riemannian map.\\

\noindent{\bf Proposition~~3.1.~}{\it Let $F_1$ be a Riemannian
submersion from a Riemannian manifold $(M_1,g_1)$ onto a Riemannian
manifold $(M_2,g_2,J)$ and  $F_2$ a slant immersion from
$(M_2,g_2,J)$ to an almost Hermitian manifold $(M_3,g_3,J)$. Then
$F_2 \circ F_1$ is a slant Riemannian map.}\\

Let $F$ be a  Riemannian map from a Riemannian manifold $(M_1,g_1)$
to an almost Hermitian manifold
$(M_2,g_2,J)$. Then for $F_*(X) \in \Gamma(range F_*)$, $X\in
\Gamma((ker F_*)^\perp)$,  we write
\begin{equation}
JF_*(X)=\phi F_*(X)+\omega F_*(X), \label{eq:3.1}
\end{equation}
where $\phi F_*(X) \in \Gamma(range F_*)$ and $ \omega F_*(X)\in
\Gamma((range F_*)^\perp)$.  Also for $V \in \Gamma((range
F_*)^\perp)$, we have
\begin{equation}
JV=\mathcal{B}V+\mathcal{C}V, \label{eq:3.2}
\end{equation}
where $\mathcal{B}V \in \Gamma(range F_*)$ and $\mathcal{C}V\in
\Gamma((range F_*)^\perp)$. We now recall from \cite{Sahin1} that
the second fundamental form $(\nabla F_*)(X,Y)$, $\forall X, Y \in
\Gamma((ker F_*)^\perp)$, of a Riemannian map has no components in
$range F_*$, i.e.
\begin{equation}
(\nabla F_*)(X,Y) \in \Gamma((range F_*)^\perp). \label{eq:3.0}
\end{equation}

 Let $F$ be a Riemannian map from
a Riemannian manifold $(M_1,g_1)$ to a
K\"{a}hler manifold $(M_2,g_2,J)$, then from (\ref{eq:3.1}),
(\ref{eq:3.2}), (\ref{eq:3.0}) and (\ref{eq:2.1}) we obtain
\begin{equation}
(\tilde{\nabla}_X \omega)F_*(Y)=\mathcal{C}(\nabla F_*)(X,Y)-(\nabla
F_*)(X,Y') \label{eq:3.3}
\end{equation}
and
\begin{equation}
F_*(\nabla^1_XY')-\phi F_*(\nabla^1_XY)=\mathcal{B}(\nabla
F_*)(X,Y)+\mathcal{S}_{\omega F_*(Y)}F_*(X)\label{eq:3.4}
\end{equation}
for $X, Y \in \Gamma((ker F_*)^\perp)$, where $(\tilde{\nabla}_X
\omega)F_*(Y)= \nabla^{^{F \perp}}_X F_*(Y)-\omega
F_*(\nabla^1_XY)$, $\nabla^1$ is the Levi-Civita connection of $M_1$
and $\phi F_*(Y)=F_*(Y')$, $Y'\in \Gamma((ker
F_*)^\perp)$.\\

  Let $F$ be a slant Riemannian map
from a Riemannian manifold $(M_1,g_1)$ to an almost Hermitian
manifold $(M_2,g_2,J)$, then we say that $\omega$ is parallel if
$(\tilde{\nabla}_X \omega)F_*(Y)=0$. We also say that $\phi$ is
parallel if $F_*(\nabla^1_XY')-\phi F_*(\nabla^1_XY)=0$ for $X,Y \in
\Gamma((ker F_*)^\perp)$ such that $\phi F_*(Y)=F_*(Y')$, $Y'\in
\Gamma((ker
F_*)^\perp)$.  \\

The proof of the following result is exactly same with slant
immersions (see \cite{CB} or \cite{Carriazo} for Sasakian case),
therefore we omit its proof.\\

\noindent{\bf Theorem~3.1.~}{\it Let $F$ be a Riemannian map from a
Riemannian manifold $(M_1,g_1)$ to  an almost Hermitian manifold
$(M_2,g_2,J)$. Then $F$ is a slant Riemannian map if and only if
there exists a constant $\lambda \in [-1,0]$ such that}
$$\phi^2 F_*(X)=\lambda F_*(X)$$
{\it for $X \in \Gamma((ker F_*)^\perp)$. If $F$ is a slant
Riemannian map, then
$\lambda=-\cos^2 \theta$.}\\

By using above theorem, it is easy to see that
\begin{eqnarray}
g_2(\phi F_*(X), \phi F_*(Y))&=&\cos^2 \theta g_1(X,Y) \label{eq:3.5}\\
g_2(\omega F_*(X), \omega F_*(Y))&=&\sin^2 \theta g_1(X,Y)
\label{eq:3.6}
\end{eqnarray}
for any $X, Y \in \Gamma((range F_*)^\perp)$.\\

We now recall the notion of adjoint map which will be useful for the
results obtained in this section (for more details see \cite{Garcia-Kupeli}). Let $F:(M_1,g_1) \longrightarrow (M_2,g_2)$
be a Riemannian map between Riemannian manifolds $(M_1,g_1)$ and
$(M_2,g_2)$. Then the adjoint map ${^*F}_*$ of $F_*$ is
characterized by $g_1(x,{^*F}_{*p_1}y)=g_2(F_{*p_1}x,y)$ for $x \in
T_{p_1}M_1$, $y \in T_{F(p_1)}M_2$ and $p_1 \in M_1$. Considering
$F^{h}_*$ at each $p_1 \in M_1$ as a linear transformation
$$F^{h}_{*{p_1}}: ((ker F_*)^{\perp}(p_1),{g_1}_{{p_1}((ker F_*)^\perp (p_1))})
\rightarrow (range F_*(p_2),{g_2}_{{p_2}(range F_*)(p_2))}),$$ we
will denote the adjoint of $F^{h}_* $ by ${^*F^{h}}_{*p_1}.$ Let
${^*F}_{*p_1}$ be the adjoint of $F_{*p_1}:(T_{p_1}M_1,{g_1}_{p_2})
\longrightarrow (T_{p_2}M_2, {g_2}_{p_2})$. Then the linear
transformation
$$({^*F}_{*p_1})^h:range F_*(p_2) \longrightarrow (ker F_*)^{\perp}(p_1)$$
defined by $({^*F}_{*p_1})^hy= {^*F}_{*p_1}y$, where $y \in
\Gamma(range F_{*p_1}), p_2=F(p_1)$, is an isomorphism and
$(F^h_{*p_1})^{-1}=({^*F}_{*p_1})^h={^*(F^h_{*p_1})}$.\\

Let $\{e_1,..,e_n\}$ be an orthonormal basis of $(ker F_*)^\perp$.
Then $\{F_*(e_1),...,F_*(e_n)\}$ is an orthonormal basis of $range
F_*$. By using (\ref{eq:3.5}) we can easily conclude that
$$\{F_*(e_1),\sec \theta \phi F_*(e_1), F_*(e_2),\sec \theta \phi e_2,...,e_p,\sec \theta \phi
e_p\}$$ is an orthonormal frame for $\Gamma(range F_*)$, where
$2p=n=rank F_*$. Then we have the following result.\\

\noindent{\bf Lemma~3.1.~}{\it  Let $F$ be a slant Riemannian map
from a Riemannian manifold $(M_1,g_1)$ to an almost Hermitian
manifold $(M_2,g_2,J)$ with $rank F_*=n$. Let
$$\{e_1,e_2,..,e_p\}, p=\frac{n}{2}$$
be a set of orthonormal vector fields in $(ker F_*)^\perp$. Then
$$\{e_1,\sec \theta {^*F}_*\phi F_*(e_1),e_2,\sec \theta {^*F}_*\phi
F_*(e_2),...,e_p,\sec \theta {^*F}_*\phi F_*(e_p)\}$$ is a local
orthonormal basis of $(ker F_*)^\perp$.}\\

\noindent{\bf Proof.~} First by direct computation we have
$$g_1(e_i,\sec \theta {^*F}_*\phi F_*(e_i))=\sec \theta
g_2(F_*(e_i),\phi F_*(e_i))=\sec \theta g_2(F_*(e_i),J
F_*(e_i))=0.$$ In a similar way, we have
$$g_1(e_i,\sec \theta {^*F}_*\phi F_*(e_j))=0.$$ Since for a Riemannian map, we have  ${^*F}_*\circ F_*= I$ (Identity map), by using
(\ref{eq:3.5}), we get
$$g_1(\sec \theta {^*F}_*\phi F_*(e_i),\sec \theta {^*F}_*\phi
F_*(e_j))=\sec^2 \theta \cos^2 \theta g_1(e_i,e_j)=\delta_{ij},$$
which gives the assertion.\\

We now denote ${^*F}_*\phi F_*$ by $Q$, then we have the following
characterization of slant Riemannian maps.\\

\noindent{\bf Corollary~3.1.~}{\it Let $F$ be a Riemannian map from
a Riemannian manifold $(M_1,g_1)$ to  an almost Hermitian manifold
$(M_2,g_2,J)$. Then $F$ is a slant Riemannian map if and only if
there exists a constant $\mu \in [-1,0]$ such that}
$$Q^2X=\mu X$$
{\it for $X \in \Gamma((ker F_*)^\perp)$. If $F$ is a slant
Riemannian map, then $\mu=-\cos^2 \theta$.}\\

\noindent{\bf Proof.~} By direct computation we have
$$Q^2X={^*F}_*\phi^2F_*(X).$$
Then proof comes from Theorem 3.1.\\

In the sequel we are going to show that the notion $\omega$ is
useful to investigate the harmonicity of slant Riemannian map. To
see this, we need the following Lemma.\\

\noindent{\bf Lemma~3.2.~}{\it Let $F$ be a slant Riemannian map
from a Riemannian manifold $(M_1,g_1)$ to  a K\"{a}hler manifold
$(M_2,g_2,J)$. If $\omega$ is parallel then we have}
\begin{equation}
(\nabla F_*)(QX,QY)=-\cos^2\,\theta\,(\nabla F_*)(X,Y)
\label{eq:3.7}
\end{equation}
{\it for $X,Y\in \Gamma((ker F_*)^\perp)$.}\\

\noindent{\bf Proof.~} If $\omega$ is parallel, then from
(\ref{eq:3.3}) we have
$$\mathcal{C}(\nabla F_*)(X,Y)=(\nabla
F_*)(X,QY).$$ Interchanging the role of $X$ and $Y$ and taking into
account that the second fundamental form is symmetric, we obtain
$$(\nabla F_*)(QX,Y)=(\nabla
F_*)(X,QY).$$ Hence we get
$$(\nabla F_*)(QX,QY)=(\nabla
F_*)(X,Q^2Y).$$ Then corollary 3.1 implies (\ref{eq:3.7}).\\

 \noindent{\bf
Theorem~3.2.~}{\it Let $F$ be a slant Riemannian map from a
Riemannian manifold $(M^{m_1}_1,g_1)$ to a K\"{a}hler manifold
$(M^{m_2}_2,g_2,J)$. If $\omega$ is parallel then $F$ is harmonic if
and
only if the distribution $(ker F_*)$ is minimal.}\\

\noindent{\bf Proof.~}  We choose an orthonormal basis of $TM_1$ as
 $$\{ v_1,...,v_{r}, e_1,\sec \theta
Qe_1,e_2,\sec \theta Qe_2,...,e_{s},\sec \theta
Qe_{s}\},\,r+2s=m_1$$ where $\{ v_1,...,v_{r}\}$ is an
orthonormal basis of $kerF_*$ and $\{ e_1, \sec \theta Qe_1$, $e_2,
\sec \theta Qe_2$,$ ...$, $e_{s},\sec \theta Qe_{s}\}$ is an
orthonormal basis of $(ker F_*)^\perp$. Since the second fundamental
form is linear in every slot, we have

$$\tau=\sum^{r}_{i=1}(\nabla F_*)(v_i,v_i)+\sum^{s}_{j=1}(\nabla
F_*)(e_j,e_j)+(\nabla F_*)(\sec\,\theta\,Qe_j,\sec\,\theta\,Qe_j).$$
Then from (\ref{eq:3.7}) we obtain
$$\tau=\sum^{r}_{i=1}(\nabla
F_*)(v_i,v_i)=-\sum^{r}_{i=1}F_*(\nabla_{v_i}v_i)$$ which proves
the assertion.\\

In the rest of this section, we investigate necessary and sufficient
conditions for a slant Riemannian map to be totally geodesic. We
recall that a differentiable map $F$ between Riemannian manifolds
$(M_1,g_1)$ and $(M_2,g_2)$ is called a totally geodesic map if
$(\nabla F_*)(X,Y)=0$ for all $X, Y \in \Gamma(TM_1)$.\\

\noindent{\bf Theorem~3.3.~}{\it Let $F$ be a slant Riemannian map
from a Riemannian manifold $(M^{m_1}_1,g_1)$ to a K\"{a}hler
manifold $(M^{m_2}_2,g_2,J)$. Then $F$ is totally geodesic if and
only if}
\begin{enumerate}
  \item [(i)]{\it the fibres are totally geodesic,}
  \item [(ii)] {\it the horizontal distribution $(ker F_*)^\perp$ is
  totally godesic,}
  \item [(iii)] {\it For $X,Y \in \Gamma((ker F_*)^\perp)$ and $V\in
  \Gamma((range F_*)^\perp)$ we have}
 $$ g_2(\mathcal{S}_{\omega
  F_*(Y)}F_*(X),\mathcal{B}V)=g_2(\nabla^{F\perp}_X\omega
  F_*(Y),\mathcal{C}V)-g_2(\nabla^{F\perp}_X\omega \phi
  F_*(Y),V).$$
\end{enumerate}

\noindent{\bf Proof.~} For $X, Y \in \Gamma((ker F_*)^\perp)$ and $V
\in \Gamma((range F_*)^\perp)$, we have
$$g_2((\nabla F_*)(X,Y),V)=g_2(\nabla^F_XF_*(Y),V).$$
Then using (\ref{eq:2.5}), (\ref{eq:3.1}) and (\ref{eq:3.2}) we
obtain
$$g_2((\nabla F_*)(X,Y),V)=-g_2(\nabla^F_XJ\phi F_*(Y),V)+g_2(\nabla^F_X\omega F_*(Y),JV) .$$
Using again (\ref{eq:3.1}), (\ref{eq:3.2}), (\ref{eq:2.2}) and
Theorem 3.1, we get
\begin{eqnarray}
g_2((\nabla
F_*)(X,Y),V)&=&\cos^2\,\theta\,g_2(\nabla^F_XF_*(Y),V)-g_2(\nabla^{F\perp}_X\omega
\phi F_*(Y),V)\nonumber\\
&&-g_2(\mathcal{S}_{\omega
F_*(Y)}F_*(X),\mathcal{B}V)\nonumber\\
&&+g_2(\nabla^{F\perp}_X\omega F_*(Y),\mathcal{C}V)\label{eq:3.9}
\end{eqnarray}
For $X, Y \in \Gamma((ker F_*)^\perp)$ and $W \in \Gamma(ker F_*)$,
from (\ref{eq:2.1}) we derive
\begin{equation}
g_2((\nabla F_*)(X,W),F_*(Y))=g_1(\nabla^1_X Y,W). \label{eq:3.10}
\end{equation}
In a similar way, for $U \in \Gamma(ker F_*)$, we obtain
\begin{equation}
g_2((\nabla F_*)(U,W),F_*(Y))=-g_1(\nabla^1_UW,Y). \label{eq:3.11}
\end{equation}
Then proof follows from (\ref{eq:3.0}), (\ref{eq:3.9}),
(\ref{eq:3.10}) and (\ref{eq:3.11}).

\section*{4.~Slant Riemannian maps and PHWC maps}
  \setcounter{equation}{0}
\renewcommand{\theequation}{4.\arabic{equation}}
In this section, we are going to show that every slant Riemannian
map is a  pseudo horizontally weakly conformal (PHWC) map, then we
investigate the conditions for slant Riemannian maps to be pseudo
horizontally homothetic map. First, we recall main definitions for
PHWC maps.  Let $F$ be a map from a Riemannian manifold $(M_1,g_1)$
to a K\"{a}hler manifold $(M_2,g_2,J)$, where $g_1$ and $g_2$ are
Riemannian metrics on $M_1$ and $M_2$, and $J$ is the complex
structure on $M_2$. For any point $p\in M_1$, we denote the adjoint
map of the tangent map $F_{*p}:T_pM_1\longrightarrow T_{F(p)}M_2$ by
$^*F_{*p}:T_{F(p)}M_2\longrightarrow T_pM_1$. If $rangeF_{*p}$ is
$J-$ invariant, then we can define an almost complex structure
$\tilde{J}_p$ on the horizontal space $(ker F_{*p})^\perp$ by
$$\tilde{J}_p=F^{-1}_{*p}\circ J_{F(p)}\circ F_{*p}.$$
If the spaces $range F_{*p}$ are $J-$ invariant for all $p$, then
the almost complex structure on $(ker F_*)^\perp$ is defined by
$$\tilde{J}=F^{-1}_{*}\circ J\circ F_{*}.$$
The map $F$ is called  pseudo-horizontally weakly conformal (PHWC)
at $p$ if and only if $range F_{*p}$ is $J$-invariant and
$g_1\mid_{(ker F_{*p})^\perp}$ is $\tilde{J}_p-$Hermitian. The map
$F$ is called pseudo-horizontally weakly conformal if and only if it
is pseudo-horizontally weakly conformal at any point of $p$. A
pseudo-horizontally weakly conformal map $F$ from a Riemannian
manifold $(M_1,g_1)$  to a K\"{a}hler manifold $(M_2,g_2,J)$ is
called pseudo-horizontally homothetic if and only if $\tilde{J}$ is
parallel in horizontal directions,i.e., $\nabla^1_X \tilde{J}=0$ for
$X\in \Gamma((ker F_*)^\perp)$ (for more details see \cite{AA}, \cite{AAB}).\\

\noindent{\bf Proposition~4.1.~}{\it Let $F$ be a slant Riemannian
map from a Riemannian manifold $(M_1,g_1)$ to an almost Hermitian
manifold $(M_2,g_2,J)$. Then $F$ is PHWC map.}\\

\noindent{\bf Proof.~} We first note that
$\tilde{J}=\sec\,\theta\,\phi$ is a complex structure on $(range
F_*)$ and $range F_*$ is invariant with respect to $\tilde{J}$. Then
we define $\hat{J}=\sec\,\theta\,Q=\sec\,\theta\,{^*F}_*\phi F_*$,
it is easy to see that $\hat{J}$ is a complex structure on $(ker
F_*)^\perp$. Thus $((ker F_*)^\perp, \hat{J})$ is an almost complex
distribution. We now consider $\hat{g}={g_1}\mid_{(ker F_*)^\perp}$,
then  by direct computation we obtain
$$\hat{g}(\hat{J}X,\hat{J}Y)=\hat{g}(X,Y)$$
for $X,Y\in \Gamma((ker F_*)^\perp)$. Thus $\hat{g}$ is $\hat{J}-$
Hermitian and $((ker F_*)^\perp, \hat{g},\hat{J})$ is almost
Hermitian distribution. Therefore $F$ is PHWC map.\\

From (\ref{eq:3.3}) and (\ref{eq:2.1}) we have the following.\\

\noindent{\bf Lemma~4.1.~}{\it Let $F$ be a slant Riemannian map
from a Riemannian manifold $(M_1,g_1)$ to a K\"{a}hler manifold
$(M_2,g_2,J)$. Then $\phi$ is parallel if and only if}
$$(\nabla F_*)(X,QY)=\nabla^{F}_X\phi
F_*(Y)-\phi F_*(\nabla^1_XY)$$ {\it for $X, Y\in \Gamma((ker
F_*)^\perp)$.}\\

We now give necessary and sufficient conditions for a slant
Riemannian map to be pseudo horizontally homothetic map.\\

\noindent{\bf Theorem~4.1.~}{\it Let $F$ be a slant Riemannian map
from a Riemannian manifold $(M_1,g_1)$ to a K\"{a}hler manifold
$(M_2,g_2,J)$. Then $F$ is pseudo horizontally homothetic map if and
only if $\phi$ is parallel and}
$$(\nabla F_*)(X,U)=0$$
{\it for $X\in \Gamma((ker F_*)^\perp)$ and $U\in \Gamma(ker F_*)$.}

\noindent{\bf Proof.~} By direct computations we have
$$(\nabla_X \hat{J})Y=\sec\, \theta\,\nabla^1_XQY-Q\nabla^1_XY.$$ Hence we obtain
$$F_*(\nabla_X \hat{J})Y=\sec\, \theta\,(F_*(\nabla^1_XQY)-\phi F_*(\nabla^1_XY)).$$ Then using (\ref{eq:2.1}) we get
\begin{equation}
F_*(\nabla_X \hat{J})Y=\sec\,\theta\,(-(\nabla
F_*)(X,QY)+\nabla^{F}_X\phi F_*(Y)-\phi F_*(\nabla^1_XY)).
\label{eq:4.1}
\end{equation}
On the other hand, since $QY$ and $U \in \Gamma(Ker F_*)$ are orthogonal, we have
$$g_1((\nabla^1_X\hat{J})Y,U)=\sec\,\theta\,g_1(QY,\nabla^1_XU).$$ Then adjoint map ${^*F}_*$ and (\ref{eq:2.1})
imply
\begin{equation}
g_1((\nabla^1_X\hat{J})Y,U)=\sec\,\theta\,g_2(\phi F_*(Y),(\nabla
F_*)(X,U)).\label{eq:4.2}
\end{equation}
Then proof comes from (\ref{eq:4.1}) and (\ref{eq:4.2}).\\

\noindent{\bf Remark~4.1.~}A harmonic map which is pseudo
horizontally weakly  is called a pseudo harmonic morphism. Pseudo
harmonic morphisms were defined in \cite{Loubeau} by Loubeau. He
also obtained that a smooth map $F$ from a Riemannian manifold $M$
to a K\"{a}hler manifold $N$ is a pseudo-harmonic morphism if and
only if $F$ pulls back local pluriharmonic functions on $N$ to local
harmonic functions on $M$. In this paper, we show that any slant
Riemannian map is a pseudo horizontally weakly conformal map
(Proposition 4.1). Also Theorem 3.2 of the present paper shows that
it is possible to obtain harmonic slant Riemannian maps under some
geometric conditions. Thus slant Riemannian maps are good candidates
for pseudo harmonic morphisms.\\

\noindent{\bf Acknowledgment.} This paper is supported by The
Scientific and Technological Council of Turkey (TUBITAK) with number
(TBAG-109T125).


\begin{thebibliography}{25}
\bibitem{Marsden} Abraham,R., Marsden, J.E., Ratiu,T., {\it Manifolds, Tensor Analysis, and Applications}, Springer-Verlag, Newyork, 1988.
\bibitem{AA}Aprodu, M.A., Aprodu, A. M., Implicitly defined harmonic PHH
submersions, Manuscripta Math. 100,(1999), 103-121.
\bibitem{AAB} Aprodu, M.A., Aprodu, A. M., Br$\hat{i}$nz$\breve{a}$nescu, A class of harmonic maps
and minimal submanifolds, Int. J. Math. 11 (2000), 1177-1191.
\bibitem{Baird-Wood} Baird, P. and  Wood, J. C. {\it Harmonic Morphisms Between Riemannian Manifolds}, London Mathematical Society Monographs, No. 29, Oxford University Press, The Clarendon Press, Oxford, 2003.
\bibitem{BBBR} Burns, D., Burstall, F. E., de Bartolomeis, P., Rawnsley, J., Stability of harmonic
maps of K\"{a}hler manifolds. J. Diff. Geom. 30, (1989), 579–594.
\bibitem{Carriazo}  Cabrerizo, J. L., Carriazo, A., Fernandez, L. M., Fernandez, M., Slant submanifolds in Sasakian manifolds. Glasg. Math. J. 42(1), (2000),125-138.
\bibitem{Chen} Chen, B. Y., Slant immersions. Bull. Austral. Math. Soc.
41 (1), (1990), 135-147.
\bibitem{CB}Chen, B. Y., {\it Geometry of slant Submanifolds}, Katholieke Universiteit Leuven, Leuven, 1990.
\bibitem{Domingo}  Chinea, D., Almost contact metric submersions. Rend. Circ. Mat. Palermo, (1985), 34(1), 89-104.
\bibitem{Falcitelli} Falcitelli, M., Ianus, S., Pastore, A. M., {\it Riemannian Submersions and Related Topics}. World Scientific, River Edge, NJ, 2004.
\bibitem{Fischer} Fischer, A. E.: Riemannian maps between
Riemannian manifolds, Contemporary math. 132, 331-366, (1992).
\bibitem{Garcia-Kupeli}Garcia-Rio, E., Kupeli,D. N.,  {\it Semi-Riemannian maps and their
Applications}, Kluwer Academic, Dortrecht, 1999.
\bibitem{Gray} Gray, A., Pseudo-Riemannian almost product manifolds and submersions, J. Math. Mech., (1967), 16, 715-737.
\bibitem{Ianus} Ianus, S., Mazzocco, R., Vilcu, G. E., Riemannian submersions from quaternionic manifolds. Acta Appl. Math., (2008), 104(1), 83-89.
\bibitem{Loubeau}Loubeau, E. Pseudo-harmonic morphisms, Internat. J. Math.
8, no. 7,(1997), 943-957.
\bibitem{LS}Loubeau, E, Slobodeanu, R., Eigenvalues of harmonic almost
submersions, Geom Dedicata  145, (2010), 103–126
\bibitem{lck}Marrero, J. C., Rocha, J., Locally conformal K\"{a}hler submersions. Geom. Dedicata, (1994), 52(3), 271-289.
\bibitem{Nore}Nore, T., Second fundamental form of a map, Ann. Mat. Pur. and Appl., 146, (1987),  281-310.
\bibitem{O'Neill} O'Neill, B., The fundamental equations of a submersion, Mich. Math. J., (1966), 13, 458-469.
\bibitem{Sahin1} \d{S}ahin, B., Invariant and anti-invariant
Riemannian maps to K\"{a}hler manifolds, Int. J. Geom. Methods Mod. Phys., vol:7, no:3 (2010), 1-19.
\bibitem{Sahin2}\d{S}ahin, B., Semi-invariant Riemannian maps to K\"{a}hler manifolds, Int. J. Geom. Methods Mod. Phys., 8(7), (2011), 1439-1454.
\bibitem{Watson}Watson, B. Almost Hermitian submersions. J. Differential Geometry, (1976), 11(1), 147-165.
\bibitem{Yano-Kon} Yano, K. and Kon, M. {\it Structures on Manifolds}, World Scientific, Singapore, 1984.

\end{thebibliography}
\end{document}